\documentclass{amsart}
\usepackage{amssymb}
\usepackage{amsfonts}

\newtheorem{theorem}{Theorem}
\theoremstyle{plain}

\newtheorem{proposition}{Proposition}
\newtheorem{remark}{Remark}

\numberwithin{equation}{section}
\input{tcilatex}

\newcommand{\dint} {\displaystyle\int}

\begin{document}
\title[Representation theoretic patterns in three dimensional Cryo-EM]{%
Representation theoretic patterns in three dimensional Cryo-Electron
Microscopy I - The intrinsic reconstitution algorithm}
\author{Ronny Hadani}
\curraddr{Department of Mathematics, University of Texas at Austin, Austin
C1200, USA. }
\author{Amit Singer}
\address{Department of Mathematics and PACM, Princeton University, Fine
Hall, Washington Road, Princeton NJ 08544-1000, USA}
\email{hadani@math.utexas.edu}
\email{amits@math.princeton.edu}
\date{October 15, 2009}

\begin{abstract}
In this paper, we describe and study a mathematical framework for
cryo-elecron microscopy. The main result, is a a proof of the admissability
(correctness) and the numerical stability of the intrinsic reconstitution
algorithm which was introduced by Singer and Shkolnisky in \cite{S}. In
addition, we explain how the various numerical observations reported in\
that work, follow from basic representation theoretic principles.
\end{abstract}

\maketitle

\section{Introduction}

"Three dimensional electron microscopy" \cite{F} is the name commonly given
to methods in which the three dimensional structure of a macromolecular
complex is obtained from the set of images taken by an electron microscope.
The most general and widespread of this methods is single-particle
reconstruction. In this method the three dimensional structure is determined
from images of randomly oriented and positioned identical macromolecular
complexes, also referred to as \textit{molecular} \textit{particles}. A
variant of this method is called cryo-electron microscopy (or cryo-EM for
short) where multitude of molecular particles are rapidly immobilized in a
thin layer of ice and maintained at liquid nitrogen \ temperature throughout
the imaging process. Single particle reconstruction from cryo-electron
microscopy images is of particular interest, since it promises to be an
entirely general technique which does not require crystallization or other
special preparation stages and is beginning \ to reach sufficient resolution
to allow the trace of polypeptide chains and the identification of residues
in protein molecules \cite{H,L,Z}.

Over the years, several methods were proposed for single particle
reconstruction from cryo-EM images. Present methods are based on the
"Angular Reconstitution" algorithm of Van Heel \cite{V}, which was also
developed independently by Vainshtein and Goncharov \cite{G}. \ However,
these methods fail with particles that are too small, cryo-EM images that
are too noisy or at resolutions where the signal-to-noise ratio becomes too
small.

\subsection{Main results}

In \cite{S}, a novel algorithm, referred to in this paper as \textit{the
intrinsic reconstitution algorithm}, for single particle reconstruction from
cryo-EM images was presented. The appealing property of this new algorithm
is that it exhibits remarkable numerical stability to noise. The
admissibility (correctness) and the numerical stability of this algorithm
were verified in an overwhelming number of numerical simulations, albeit, a
formal justification was still missing. In this paper, we prove the
admissibility and the numerical stability of the intrinsic reconstitution
algorithm. The proof relies on the study of a certain operator $C$, of
geometric origin, referred to as the \textit{common lines operator}.
Specifically,

\begin{itemize}
\item Admissibility, depends, among other things, on the fact that the
maximal eigenspace of $C$ is three dimensional.

\item Numerical stability, depends on the existence of a spectral gap which
separates the maximal eigenvalue $\lambda _{\max }$ from the rest of the
spectrum.
\end{itemize}

In this regard, the main technical result of this paper is a complete
description of the spectral properties of the common lines operator. In the
course, we describe a formal mathematical framework for cryo-EM which
explains how the various numerical observations reported in \cite{S} follow
from basic representation theoretic principles, thus putting that work on
firm mathematical grounds.

The remainder of the introduction is devoted to a detailed description of
the intrinsic reconstitution algorithm and to the explanation of the main
ideas and results of this paper.

\subsection{Mathematical model}

When modeling the mathematics of cryo-EM its more convenient to think of a
fixed macromolecular complex which is observed from different directions by
the electron microscope. In more details, the macromolecular complex is
modeled by a real valued function $\phi $ on a three dimensional Euclidian
vector space $V\simeq
\mathbb{R}
^{3}$, which \ describes the electric potential due to the charge density in
the complex. A viewing direction of the electron microscope is modeled by a
point $x$ on the unit sphere $X=S\left( V\right) $. The interaction of the
beam from the electron microscope with the complex is modeled by a real
valued function $R_{x}$ on the orthogonal plane $\mathrm{P}_{x}=x^{\bot }$,
given by the Radon transform of $\phi $ along the direction $x$, that is
\begin{equation*}
R_{x}\left( v\right) =\int \limits_{L_{x}}\phi \left( v+l\right) dl,
\end{equation*}%
for every $v\in \mathrm{P}_{x}$, where $L_{x}$ is the line passing thorough $%
x$ and $dl$ is the Euclidian measure on $L_{x}$.

The data collected from the experiment is a collection of Radon images $%
R_{x}:\mathrm{P}_{x}\rightarrow
\mathbb{R}
$, $x\in X_{N}$; where $X_{N}\subset X$ consists of $N$ points. The main
empirical assumption is

\medskip

\textbf{Empirical assumption: }The points $x\in X_{N}$ are distributed
independently and uniformly at random.

\medskip

We emphasize that, in practice, the embeddings $i_{x}:\mathrm{P}%
_{x}\hookrightarrow V$, $x\in X_{N}$ are not known. What one is given is the
collection of Radon images $\left \{ R_{x}:\mathrm{P}_{x}\rightarrow
\mathbb{R}
:x\in X_{N}\right \} $, where each plane $\mathrm{P}_{x},x\in X_{N}$ should
be considered as an abstract Euclidian plane.

The main problem of cryo-EM is to reconstruct the orthogonal embeddings $%
\left \{ i_{x}:\mathrm{P}_{x}\hookrightarrow V:x\in X_{N}\right \} $ from
the Radon images $\left \{ R_{x}:\mathrm{P}_{x}\rightarrow
\mathbb{R}
:x\in X_{N}\right \} $. We will refer to this problem as the \textit{cryo-EM
reconstruction problem} and note that, granting its solution, the function $%
\phi $ can be computed (approximately) using the inverse Radon transform.

\subsection{The Fourier slicing property and the common lines datum}

The first step of the reconstruction is to extract from the Radon images a
certain linear algebra datum, referred to as \textit{the common lines datum}%
, which captures a basic relation in three dimensional Euclidian geometry.

First one notes that for every $x\in X_{N}$, the following relation holds
\begin{equation}
\widehat{R}_{x}=\widehat{\phi }_{|\mathrm{P}_{x}},
\label{Fourier_slicing_eq}
\end{equation}%
where the operation $\widehat{\left( -\right) }$ on the left hand side
denotes the Euclidian Fourier transform on the plane $\mathrm{P}_{x}$ and
the operation $\widehat{\left( -\right) }$ on the right hand side denotes
the Euclidian Fourier transform on $V$. This relation follows easily from
the standard properties of the Fourier transform and is sometimes referred
to as the \textit{Fourier slicing property }(see \cite{N}). The key
observation, first made by Klug (see \cite{K}), is that (\ref%
{Fourier_slicing_eq}) implies, for every pair of different points $x,y\in
X_{N}$, that the functions $\widehat{R}_{x}$ and $\widehat{R}_{y}$ must
agree on the line of intersection (common line), that is
\begin{equation*}
\widehat{R}_{x}{}_{|\mathrm{P}_{x}\cap \mathrm{P}_{y}}=\widehat{R}_{y}{}_{|%
\mathrm{P}_{x}\cap \mathrm{P}_{y}}.
\end{equation*}

Hence, if the function $\phi $ is generic enough then one can compute from
each pair of images $\widehat{R}_{x}$ and $\widehat{R}_{y}$ the operator $%
C_{N}\left( x,y\right) :\mathrm{P}_{y}\rightarrow \mathrm{P}_{x}$, which
identifies the line of intersection between the two planes. Formally, this
operator is given by the composition $C_{x,y}\circ C_{y,x}^{T}$, where $%
C_{x,y}$ and $C_{y,x}$ are the tautological embeddings
\begin{eqnarray*}
C_{x,y} &:&\mathrm{P}_{x}\cap \mathrm{P}_{y}\hookrightarrow \mathrm{P}_{x},
\\
C_{y,x} &:&\mathrm{P}_{x}\cap \mathrm{P}_{y}\hookrightarrow \mathrm{P}_{y}.
\end{eqnarray*}

\subsection{The intrinsic reconstitution algorithm}

The intrinsic reconstitution algorithm reconstructs the orthogonal
embeddings $\left \{ i_{x}:\mathrm{P}_{x}\rightarrow V:x\in X_{N}\right \} $
from the common lines datum $\left \{ C_{N}\left( x,y\right) :\left(
x,y\right) \in X_{N}\times X_{N}\right \} $. The crucial step is to
construct an \textbf{intrinsic} model of the three dimensional Euclidian
vector space $V$ which is expressed solely in terms of the common lines
datum. The construction proceeds in four steps:

\begin{description}
\item[Ambient vector space] We define the $2N$ dimensional Euclidian vector
space%
\begin{equation*}
\mathcal{H}_{N}=\bigoplus \limits_{x\in X_{N}}\mathrm{P}_{x}\text{.}
\end{equation*}

\item[Common lines operator] We define the symmetric operator $C_{N}:%
\mathcal{H}_{N}\rightarrow \mathcal{H}_{N}$ by%
\begin{equation*}
C_{N}\left( s\right) \left( x\right) =\frac{1}{\left \vert X_{N}\right \vert
}\sum \limits_{y\in X_{N}}C_{N}\left( x,y\right) \left( s\left( y\right)
\right) .
\end{equation*}

\item[Intrinsic model] We define the Euclidian subspace\footnote{%
The condition $\lambda >1/3$ in the definition of $\mathbb{V}_{N}$ will be
clarified when we will discuss the spectral gap property in the next
subsection.}
\begin{equation*}
\mathbb{V}_{N}=\bigoplus \limits_{\lambda >1/3}\mathcal{H}_{N}\left( \lambda
\right) ,
\end{equation*}%
where $\mathcal{H}_{N}\left( \lambda \right) $ denote the eigenspace of $%
C_{N}$ associated with the eigenvalue $\lambda $.

\item[Intrinsic maps] For every $x\in X_{N}$, we define the map
\begin{equation*}
\varphi _{x}=\sqrt{2/3}\cdot \left( \mathrm{Pr}_{x}\right) ^{t}:\mathrm{P}%
_{x}\rightarrow \mathbb{V}_{N},
\end{equation*}%
where $\mathrm{Pr}_{x}:\mathbb{V}_{N}\rightarrow \mathrm{P}_{x}$ is the
orthogonal projection on the component $\mathrm{P}_{x}$.
\end{description}

\medskip

The fact that the vector space $\mathbb{V}_{N}$ is of the right dimension is
granted by the following theorem:

\begin{theorem}
\label{dim_thm}For sufficiently large $N$, we have
\begin{equation*}
\dim \mathbb{V}_{N}=3.
\end{equation*}
\end{theorem}

The fact that the collection of maps $\left \{ \varphi _{x}:\mathrm{P}%
_{x}\rightarrow \mathbb{V}_{N}:x\in X_{N}\right \} $ solves the cryo-EM
reconstruction problem is the content of the following theorem:

\begin{theorem}
\label{sol_thm}There exists an (approximated) isometry $\tau _{N}:V\approx
\mathbb{V}_{N}$ which satisfies the following property:
\begin{equation*}
\tau _{N}\circ i_{x}=\varphi _{x},
\end{equation*}%
for every $x\in X_{N}$.
\end{theorem}

\bigskip

\begin{remark}
Theorem \ref{sol_thm} implies that the vector space $V$ equipped with the
tautological embeddings $\left \{ i_{x}:x\in X_{N}\right \} $ is isomorphic
to the intrinsic vector space $\mathbb{V}_{N}$ equipped with the mappings $%
\left \{ \varphi _{x}:x\in X_{N}\right \} $. Hence, for all practical
purposes they are indistinguishable. This is, to our judgement, an elegant
formal example, realizing the general phylosophy about the appearence of
strucuture form data.
\end{remark}

\subsection{Analytic set-up}

The proofs of Theorems \ref{dim_thm} and \ref{sol_thm} are based on an
approximation argument of the discrete from the continuos, which we are
going to explain next.

Let $\mathfrak{H}\rightarrow X$ be the vector bundle on the unit sphere
whose fiber at a point $x\in X$ is the plane $\mathrm{P}_{x}=x^{\bot }$ and
let $\mathcal{H}$ denote the vector space of smooth global sections $\Gamma
\left( X,\mathfrak{H}\right) $. The vector space $\mathcal{H}$ is equipped
with a inner product, given by
\begin{equation*}
\left( s_{1},s_{2}\right) =\int \limits_{x\in X}B\left( s_{1}\left( x\right)
,s_{2}\left( x\right) \right) dx,
\end{equation*}%
where $dx$ is the Haar measure on the unit sphere. The common lines datum
can be used to form a kernel of a symmetric integral operator $C:\mathcal{%
H\rightarrow H}$ which is given by
\begin{equation*}
C\left( s\right) \left( x\right) =\int \limits_{y\in X}C\left( x,y\right)
\left( s\left( y\right) \right) dy,
\end{equation*}

The main difference from the discrete scenario is that, here, the space $%
\mathcal{H}$ supports a representation of the orthogonal group $O\left(
V\right) $ and the main observation is that the operator $C$ commutes with
the group action. This enables to understand the operator $C$ in terms of
the representation theory of the orthogonal group and, consequently, to
compute its spectrum and to describe the associated eigenspaces. In this
regard, the main technical result of this paper is

\begin{theorem}
\label{spec_thm}The operator $C$ admits a discrete spectrum $\lambda
_{1},\lambda _{2},...,\lambda _{n},..,\lambda _{\infty }=0$ such that%
\begin{equation*}
\lambda _{n}=\frac{\left( -1\right) ^{n-1}}{n\left( n+1\right) }.
\end{equation*}

Moreover, $\dim \mathcal{H}\left( \lambda _{n}\right) =2n+1$.
\end{theorem}

An immediate consequence of Theorem \ref{spec_thm} is that the maximal
eigenvalue of $C$ is $\lambda _{\max }=1/2$, its multiplicity is equal $3$
and there exists a \textbf{spectral gap} of $\lambda _{1}-\lambda _{3}=5/12$%
, which separates it from the rest of the spectrum. Now, consider the vector
space $\mathbb{V}=\mathcal{H}\left( \lambda _{\max }\right) $ and define the
maps
\begin{equation*}
\varphi _{x}=\sqrt{2/3}\cdot \left( ev_{x}\right) ^{t}:\mathrm{P}%
_{x}\rightarrow \mathbb{V}.
\end{equation*}

The second main result of this paper asserts that the vector space $V$
equipped with the collection of tautological embeddings $\left \{ i_{x}:x\in
X\right \} $ is isomorphic to the vector space $\mathbb{V}$ equipped with
the mappings $\left \{ \varphi _{x}:x\in X\right \} $ and, in addition, this
isomorphism is proportional to the canonical morphism $\alpha
_{can}:V\rightarrow \mathcal{H}$ which sends a vector $v\in V$ to the
section $\alpha _{can}\left( v\right) \in \mathcal{H}$ whose value at the
point $x$ is the orthogonal projection of $v$ on the plane $\mathrm{P}_{x}$.
All of this is summarized in the following theorem:

\begin{theorem}
\label{sol_analytic_thm}The morphism $\tau =\sqrt{3/2}\cdot \alpha _{can}$
maps the vector space $V$ isometrically onto $\mathbb{V}\subset \mathcal{H}$%
. Moreover,
\begin{equation*}
\tau \circ i_{x}=\varphi _{x},
\end{equation*}%
for every $x\in X$.
\end{theorem}

\subsubsection{Proof of Theorems \protect \ref{dim_thm} and \protect \ref%
{sol_thm}}

The proof is based on the following three assertions:

\begin{description}
\item[Assertion 1] The vector space $\mathcal{H}_{N}$ approximates the
vector space $\mathcal{H}$.

\item[Assertion 2] The operator $C_{N}$ approximates the integral operator $%
C $.

\item[Assertion 3] The vector space $\mathbb{V}_{N}$ approximates the vector
space $\mathbb{V}$.
\end{description}

The validity of the first two assertions depends on our principal assumption
that the points in $X_{N}$ are chosen independently and uniformly at random,
which, implies that the Haar measure on $X$ is approximated by the
(normalized) counting measure on $X_{N}$. The validity of the third
assertion depends also on the existence of a spectral gap for the operator $%
C $ which implies that the maximal eigenspace can be computed in a
numerically stable manner. Consequently, Theorems \ref{dim_thm} and \ref%
{sol_thm} follow from Theorem \ref{sol_analytic_thm}.

\begin{remark}
The reconstruction of the orthogonal maps $\left \{ i_{x}:x\in X_{N}\right
\} $ is a non-linear problem because of the orthogonality constraint. One of
the appealing properties of the intrinsic reconstitution algorithm is that
it reduces this problem to a problem in linear algebra - the computation of
the maximal eigenspace of a linear operator. Another appealing property is
its remarkable stability to noise which, using a bit of random matrix theory
arguments (see \cite{S}), follows \ from the spectral gap property. Other
existing reconstruction methods, like the angular reconstitution method (see
\cite{V} and \cite{G}), do not enjoy this important stability property.
\end{remark}

\subsection{Structure of the paper}

The paper consists of three sections except of the introduction.

\begin{itemize}
\item In Section \ref{prelim_sec}, we begin by introducing the basic
analytic set-up which underlies cryo-EM. Then, we proceed to formulate the
main results of this paper, which are: A complete description of the
spectral properties of the common lines operator $C$ (Theorem \ref%
{spectral_thm}) and the admissibility of the intrinsic reconstitution
algorithm (Theorem \ref{char_thm}).

\item In Section \ref{spectral_sec}, we prove Theorem \ref{spectral_thm}, in
particular, we develop all the representation theoretic machinery which is
needed for the proof.

\item Finally, in Appendix \ref{proofs_sec}, we give the proofs of all
technical statements which appeared in the previous sections.
\end{itemize}

\bigskip

\textbf{Acknowledgement: }\textit{The first author would like to thank
Joseph Bernstein for many helpful discussions concerning the mathematical
aspects of this work. Also, he would like to thank the MPI institute at Bonn
where several parts of this work were concluded during the summer of 2009.
We thank Shamgar Gurevich for carefully reading the manuscript and giving
various corrections and remarks. The second author is partially supported by
Award Number R01GM090200 from the National Institute of General Medical
Sciences. The content is solely the responsibility of the authors and does
not necessarily represent the official views of the National Institute of
General Medical Sciences or the National Institutes of Health.}

\section{Preliminaries and main results\label{prelim_sec}}

\subsection{Set up}

Let $\left( V,B\right) $ be a three dimensional Euclidian vector space over $%
\mathbb{R}
$. The reader can take $V$ $=$ $%
\mathbb{R}
^{3}$ equipped with the standard inner product $B_{std}:%
\mathbb{R}
^{3}\times
\mathbb{R}
^{3}\rightarrow
\mathbb{R}
$. Let $O\left( V\right) =O\left( V,B\right) $ denote the group of
orthogonal transformations with respect to the inner product $B$; let $%
SO\left( V\right) \subset O\left( V\right) $ denote the subgroup of
orthogonal transformation which have determinant one; let $\theta \in
O\left( V\right) $ denote the element $-Id$. Let $S\left( V\right) $ denote
the unit sphere in $V$, that is, $S\left( V\right) =\left \{ v\in V:B\left(
v,v\right) =1\right \} $.

\subsection{The vector bundle of planes}

Let $\mathfrak{H}\rightarrow S(V)$ be the real vector bundle with fibers $%
\mathfrak{H}_{|x}=x^{\bot }$ and let $\mathcal{H=}\Gamma \left( S\left(
V\right) ,\mathfrak{H}\right) $ denote the space of smooth global sections.
The vector bundle $\mathfrak{H}$ admits a fiberwise Euclidian structure
induced from the one on $V$, which in turns yields a (pre) Euclidian
structure on $\mathcal{H}$ (here, the word pre just means that $\mathcal{H}$
is not complete). In general, in this paper we will not distinguish between
an Euclidian/Hermitian vector space and its completion and the correct
choice between the two will be clear from the context.

In addition, $\mathfrak{H}$ admits a natural $O\left( V\right) $ equivariant
structure which induces an orthogonal action of $O\left( V\right) $ on the
space of global sections $\mathcal{H}$ which sends a section $s\in \mathcal{H%
}$ to a section $g\cdot s$, given by
\begin{equation*}
\left( g\cdot s\right) \left( x\right) =gs\left( g^{-1}x\right) ,
\end{equation*}%
for every $x\in S(V)$. This makes $\mathcal{H}$ into an Euclidian
representation of $O\left( V\right) $. We will also consider the
complexified vector bundle $%
\mathbb{C}
\mathfrak{H}$ and its space of global sections $%
\mathbb{C}
\mathcal{H=}\Gamma \left( S\left( V\right) ,%
\mathbb{C}
\mathfrak{H}\right) $. The vector bundle $%
\mathbb{C}
\mathfrak{H}$ is equipped with an Hermitian inner product induced from the
Hermitian product $\left \langle -,-\right \rangle $ on $%
\mathbb{C}
V$ which is given by
\begin{equation*}
\left \langle u,v\right \rangle =B\left( \overline{u},v\right) \text{,}
\end{equation*}%
where $\overline{\left( -\right) }:%
\mathbb{C}
V\rightarrow
\mathbb{C}
V$ is the Galois conjugation. Consequently, $%
\mathbb{C}
\mathcal{H}$ is a (pre) Hemitian vector space supporting a (real) unitary
representation of the group $O\left( V\right) $.

\subsection{The operator of common lines}

We define an integral operator $C:\mathcal{H}\rightarrow \mathcal{H}$
capturing a basic relation in three dimensional Euclidian geometry.

The operator $C$ is defined as follows: For every pair of points $x,y\in
S\left( V\right) $ such that $x\neq \pm y$, consider the intersection $%
x^{\bot }\cap y^{\bot }$ of the corresponding orthogonal planes, which is a
line in $V$. There are two tautological embeddings
\begin{eqnarray*}
C_{x,y} &:&x^{\bot }\cap y^{\bot }\hookrightarrow x^{\bot }, \\
C_{y,x} &:&x^{\bot }\cap y^{\bot }\hookrightarrow y^{\bot }\text{.}
\end{eqnarray*}

Using these embeddings we can define a rank one operator $C\left( x,y\right)
:\mathfrak{H}_{|y}\rightarrow \mathfrak{H}_{|x}$, given by the composition $%
C_{x,y}\circ C_{y,x}^{t}$. The collection $\left \{ C\left( x,y\right)
:x\neq \pm y\right \} $ yields a well defined smooth section of $\mathfrak{H}%
\boxtimes \mathfrak{H}^{\ast }$ on the complement of the union of the
diagonal and the anti-diagonal. It is not difficult to verify that this
section extends to a distribution section $C\in \Gamma ^{\prime }\left(
S\left( V\right) \times S\left( V\right) ,\mathfrak{H}\boxtimes \mathfrak{H}%
^{\ast }\right) $ which, in turns, establishes a kernel for an integral
operator $C:\mathcal{H}\rightarrow \mathcal{H}$ given by
\begin{equation*}
C\left( s\right) \left( x\right) =\dint \limits_{y\in S\left( V\right)
}C\left( x,y\right) \left( s\left( y\right) \right) dy,
\end{equation*}%
for every $s\in \mathcal{H}$, where we take $dy$ to be the normalized Haar
measure on the sphere. Since $C\left( x,y\right) =C\left( y,x\right) ^{t}$,
this implies that $C$ is a symmetric operator. In addition, it is evident
that $C$ commutes with the $O\left( V\right) $ action, namely $C\left(
g\cdot s\right) =g\cdot C\left( s\right) $ for every $s\in \mathcal{H}$ and $%
g\in O\left( V\right) $.

The operator $C$ is referred to as the \textit{operator of} \textit{common
lines} and the main technical part of this paper will be devoted to the
investigation of this operator.

\subsection{Main results}

Our goal is to describe an\textbf{\ intrinsic} model of the Euclidian vector
space $V$ which is expressed in the terms of the common lines operator and
the Euclidian structure of $\mathcal{H}$ alone.

The main technical result of this paper is

\begin{theorem}
\label{spectral_thm}The operator $C$ admits a discrete spectrum $\lambda
_{1},\lambda _{2},..,\lambda _{n},..\lambda _{\infty }=0$, such that%
\begin{equation*}
\lambda _{n}=\frac{\left( -1\right) ^{n-1}}{n\left( n+1\right) }.
\end{equation*}

Moreover, $\dim \mathcal{H}\left( \lambda _{n}\right) =2n+1$.
\end{theorem}

For a proof, see Section \ref{spectral_sec}.

\medskip

\textbf{Intrinsic model: }Take $\mathbb{V}=\mathcal{H}\left( \lambda _{\max
}\right) $ to be the maximal eigenspace of $C$.

\medskip

There are two immediate implications of Theorem \ref{spectral_thm} that we
will consider: \medskip \

\begin{itemize}
\item The vector space $\mathbb{V}$ is three dimensional.

\item There exists a spectral gap of $\lambda _{1}-\lambda _{3}=5/12$ which
separates $\lambda _{\max }$ from the rest of the spectrum.
\end{itemize}

Fix $r=\sqrt{3/2}$. Let $ev_{x}:\mathcal{H}\rightarrow x^{\bot }$ denote the
evaluation morphism at the point $x\in S\left( V\right) $. For every $x\in
S\left( V\right) $, define the morphism
\begin{equation*}
\varphi _{x}=r^{-1}\cdot \left( ev_{x|\mathbb{V}}\right) ^{t}:x^{\bot
}\rightarrow \mathbb{V}\text{.}
\end{equation*}

Let $\alpha _{can}:V\rightarrow \mathcal{H}$ be the canonical morphism which
sends a vector $v\in V$ to the section $\alpha _{can}\left( v\right) \in
\mathcal{H}$, defined by%
\begin{equation*}
\alpha _{can}\left( v\right) \left( x\right) =\mathrm{Pr}_{x}\left( v\right)
,
\end{equation*}%
where $\mathrm{Pr}_{x}$ is the orthogonal projection on the plane $x^{\bot }$%
. The morphism $\alpha _{can}$ is a morphism of representations of the
orthogonal group $O\left( V\right) $.

Finally, define%
\begin{equation*}
\tau =r\cdot \alpha _{can}:V\rightarrow \mathcal{H}\text{.}
\end{equation*}

The morphism $\tau $ identifies the Euclidian vector space $V$ together with
the tautological embeddings $\left \{ i_{x}:x^{\bot }\rightarrow V\right \} $
with the Euclidian vector space $\mathbb{V}$ together with the maps $\left
\{ \varphi _{x}:x^{\bot }\rightarrow \mathbb{V}\right \} $. All of this is
summarized in the following theorem:

\begin{theorem}
\label{char_thm}The map $\tau $ maps $V$ isometrically onto $\mathbb{V}$.
Moreover,
\begin{equation*}
\tau \circ i_{x}=\varphi _{x},
\end{equation*}%
for every $x\in S\left( V\right) $.
\end{theorem}

For a proof, see Appendix \ref{proofs_sec} (the proof uses the results and
terminology of Section \ref{spectral_sec}).

\section{Spectral analysis of the operator of common lines\label%
{spectral_sec}}

\subsection{Set-up}

It will be convenient to extend the set-up a bit.

\subsubsection{Auxiliary vector bundles}

We introduce the following auxiliary vector bundles on $S\left( V\right) $.
Let $\mathfrak{N}\rightarrow S(V)$ be the vector bundle of normal lines with
fibers $\mathfrak{N}_{|x}=%
\mathbb{R}
x$ and let $\mathcal{N}=\Gamma \left( S\left( V\right) ,\mathfrak{N}\right) $
denote the corresponding space of global sections. Let $V_{S\left( V\right)
} $ denote the trivial vector bundle with fiber at each point equal $V$ and
let $\mathcal{V}$ $=\Gamma \left( S\left( V\right) ,V_{S\left( V\right)
}\right) =\mathcal{F\otimes }V$ where $\mathcal{F}=C^{\infty }\left( S\left(
V\right) ,%
\mathbb{R}
\right) $, we have
\begin{equation*}
\mathcal{V}=\mathcal{H}\oplus \mathcal{N}\text{.}
\end{equation*}

All the vector bundles are equipped with a fiberwise Euclidian structure
which is induced from the one on $V$ and consequently the spaces of global
sections are Euclidian. In addition, all the vector bundles are equipped
with a natural $O\left( V\right) $ equivariant structure which is compatible
with the Euclidian structure and consequently the spaces of global sections
form Euclidian representations of the group $O\left( V\right) $. We will
consider these spaces as representations of the subgroup $SO\left( V\right) $
and remember also the action of the special element $\theta \in O\left(
V\right) $ which commutes with the action of $SO\left( V\right) $.

\subsubsection{The operator of orthographic lines}

We define an integral operator $O:\mathcal{H}\rightarrow \mathcal{H}$ which
we refer to as the\textit{\ operator of} \textit{orthographic lines. }This
operator captures another basic relation in three dimensional Euclidian
geometry which, in some sense, stands in duality with the common lines
relation.

The operator $O$ is defined by the following kernel: For every pair of
points $x,y\in S\left( V\right) $ such that $x\neq \pm y$, consider the pair
of unit vectors
\begin{eqnarray*}
o_{x,y} &=&\frac{\mathrm{Pr}_{x}\left( y\right) }{\left \Vert \mathrm{Pr}%
_{x}\left( y\right) \right \Vert }\in x^{\bot }, \\
o_{y,x} &=&\frac{\mathrm{Pr}_{y}\left( x\right) }{\left \Vert \mathrm{Pr}%
_{y}\left( x\right) \right \Vert }\in y^{\bot }\text{.}
\end{eqnarray*}

In words, the vector $o_{x,y}$ is the normalized projection of the vector $y$
on the plane $x^{\bot }$ and similarly the vector $o_{y,x}$ is the
normalized projection of the vector $x$ on the plane $y^{\bot }$. Define a
rank one operator $O\left( x,y\right) :y^{\bot }\rightarrow x^{\bot }$ by
\begin{equation*}
O\left( x,y\right) \left( v\right) =B\left( o_{y,x},v\right) o_{x,y},
\end{equation*}%
for every $v\in y^{\bot }$. The collection $\left \{ O\left( x,y\right)
:x\neq \pm y\right \} $ yields a well defined smooth section of $\mathfrak{H}%
\boxtimes \mathfrak{H}^{\ast }$ on the complement of the union of the
diagonal and the anti-diagonal which extends to a distribution section $O\in
\Gamma ^{\prime }\left( S\left( V\right) \times S\left( V\right) ,\mathfrak{H%
}\boxtimes \mathfrak{H}^{\ast }\right) $ which, in turns, yields a symmetric
integral operator $O:\mathcal{H}\rightarrow \mathcal{H}$ which commutes with
the $O\left( V\right) $ action.

Given a pair of unit vectors $x\neq \pm y$, the following observations are
in order.

\begin{itemize}
\item The orthographic lines $%
\mathbb{R}
o_{x,y}\in x^{\bot }$ and $%
\mathbb{R}
o_{y,x}\in y^{\bot }$ are orthogonal to the common line $x^{\bot }\cap
y^{\bot }$.

\item The kernel $O\left( x,y\right) $ satisfy $O\left( x,-y\right) =O\left(
-x,y\right) =-1\cdot O\left( x,y\right) $ which means that $O\left(
x,y\right) $ depends on the choice of the unit vectors $x,y$ and not only on
the planes $x^{\bot },y^{\bot }$. This should be contrasted with the
analogue property of the kernel $C\left( x,y\right) $ which satisfies $%
C\left( x,-y\right) =C\left( -x,y\right) =C\left( x,y\right) $.
\end{itemize}

\subsubsection{The operator of parallel translation}

We define the integral operator $T=C-O:\mathcal{H}\rightarrow \mathcal{H}$.
\ For every $x,y\in S\left( V\right) $, such that $x\neq \pm y$, the kernel $%
T\left( x,y\right) :y^{\bot }\rightarrow x^{\bot }$ is a full rank operator
and it is not difficult to verify that $T\left( x,y\right) $ is the operator
of parallel translation along the unique geodesic (large circle) connecting
the point $y$ with $x$. Consequently, we will refer to $T$ as the operator
of \textit{parallel translations}.

The strategy that we are going to follow is to study the spectral properties
of the operator $T$, from which, as it turns out, the spectral properties of
the operators $C$ and $O$ can be derived.

\subsection{Isotypic decompositions}

The spaces $\mathcal{H},\mathcal{V},\mathcal{N}$ and $\mathcal{F}$ form
Euclidian representations of the group $SO\left( V\right) $ and as such
decompose into isotypic components\footnote{%
We remind the reader that an isotypic component is a representation which is
a direct sum of copies a single irreducible representation.
\par
{}}%
\begin{eqnarray*}
\mathcal{H} &\mathcal{=}&\bigoplus \limits_{n=0}^{\infty }\mathcal{H}_{n}, \\
\mathcal{N} &\mathcal{=}&\bigoplus \limits_{n=0}^{\infty }\mathcal{N}_{n}, \\
\mathcal{V} &\mathcal{=}&\bigoplus \limits_{n=0}^{\infty }\mathcal{V}_{n}, \\
\mathcal{F} &\mathcal{=}&\bigoplus \limits_{n=0}^{\infty }\mathcal{F}_{n}%
\text{,}
\end{eqnarray*}%
where we use the subscript $n$ to denote the isotypic component which
consists of the unique irreducible representation of $SO\left( V\right) $ of
dimension $2n+1$. In addition, the element $\theta \in O\left( V\right) $
acts on all these spaces, thus decompose them into a direct sum of two
components
\begin{eqnarray*}
\mathcal{H} &=&\mathcal{H}^{+}\oplus \mathcal{H}^{-}, \\
\mathcal{N} &=&\mathcal{N}^{+}\oplus \mathcal{N}^{-}, \\
\mathcal{V} &=&\mathcal{V}^{+}\oplus \mathcal{V}^{-}, \\
\mathcal{F} &=&\mathcal{F}^{+}\oplus \mathcal{F}^{-},
\end{eqnarray*}%
where we use the superscript $+$ to denote the component on which $\theta $
acts as $Id$ and the superscript $-$ to denote the component on which $%
\theta $ acts as $-Id$. We will refer to the $+$ component as the \underline{%
symmetric} component and to the $-$ component as the \underline{%
anti-symmetric} component.

The following theorem summarizes the properties of these decompositions
which we will require in the sequel.

\begin{theorem}
\label{dec_thm}The following properties hold

\begin{enumerate}
\item Each isotypic component $\mathcal{F}_{n}$ is an irreducible
representation. Moreover, $\mathcal{F}_{n}=\mathcal{F}_{n}^{+}$ when $n$ is
even and $\mathcal{F}_{n}=\mathcal{F}_{n}^{-}$ when $n$ is odd.

\item Each isotypic component $\mathcal{N}_{n}$ is an irreducible
representation. Moreover, $\mathcal{N}_{n}=\mathcal{N}_{n}^{+}$ when $n$ is
even and $\mathcal{N}_{n}=\mathcal{N}_{n}^{-}$ when $n$ is odd.

\item The isotypic component $\mathcal{H}_{0}=0$ and each isotypic component
$\mathcal{H}_{n}$, $n\geq 1$ decomposes under $\theta $ into a direct sum of
two irreducible representations $\mathcal{H}_{n}^{+}\oplus \mathcal{H}%
_{n}^{-}$.

\item The isotypic component $\mathcal{V}_{0}$ is equal to the symmetric
trivial representation $\mathbf{1}^{+}$ and each isotypic component $%
\mathcal{V}_{n}$, $n\geq 1$ decomposes under $\theta $ into a direct sum of
three irreducible representations $\mathcal{H}_{n}^{+}\oplus \mathcal{H}%
_{n}^{-}\oplus \mathcal{N}_{n}^{?}$ where
\begin{equation*}
?=\left \{
\begin{array}{cc}
+ & n\text{ even} \\
- & n\text{ odd}%
\end{array}%
\right. .
\end{equation*}
\end{enumerate}
\end{theorem}

For a proof, see Appendix \ref{proofs_sec}.

Since the operators $C,O$ and $T$ commute with the $O\left( V\right) $
action, they preserve all the above decompositions.

\begin{proposition}
\label{dec_prop}The following properties hold:

\begin{itemize}
\item The operator $T$ acts as scalar operator on $\mathcal{H}_{n}$,
moreover, $T_{|\mathcal{H}_{n}}=\lambda _{n}Id$ where $\lambda _{n}\neq 0$.

\item The isotypic component $\mathcal{H}_{n}^{+}\subset \ker C$, moreover, $%
C_{|\mathcal{H}_{n}^{-}}=\lambda _{n}Id$ .

\item The isotypic component $\mathcal{H}_{n}^{-}\subset \ker O$, moreover, $%
O_{|\mathcal{H}_{n}^{+}}=-\lambda _{n}Id$.
\end{itemize}
\end{proposition}

For a proof, see Appendix \ref{proofs_sec}$.$

The rest of this section is devoted to the computation of the eigenvalues $%
\lambda _{n}$.

\subsection{Computation of the eigenvalues}

Fix a point $x\in S\left( V\right) $. Let $T_{x}=\left \{ g\in SO\left(
V\right) :gx=x\right \} $ be the subgroup of rotations around $x$. Let $%
\left( H,E,F\right) \in $ $%
\mathbb{C}
Lie\left( SO\left( V\right) \right) $ be an an $sl_{2}$ triple associated
with $T_{x}$.

\subsubsection{Spherical decomposition}

For each $n\geq 1$, the complexified (Hilbert) space $%
\mathbb{C}
\mathcal{H}_{n}$ admits an isotypic decomposition with respect to the action
of $T_{x}$%
\begin{equation*}
\mathbb{C}
\mathcal{H}_{n}=\bigoplus \limits_{m=-n}^{n}\mathcal{H}_{n}^{m},
\end{equation*}%
where $H$ acts on $\mathcal{H}_{n}^{m}$ by $2mId$. Since $%
\mathbb{C}
\mathcal{H}_{n}=%
\mathbb{C}
\mathcal{H}_{n}^{+}\oplus
\mathbb{C}
\mathcal{H}_{n}^{-}$, each $\mathcal{H}_{n}^{m}$ is two dimensional.

\subsubsection{Strategy}

Given a section $u_{n}\in
\mathbb{C}
\mathcal{H}_{n}$, by Proposition \ref{dec_prop}, $Tu_{n}=\lambda _{n}u_{n}$.
If, in addition, $u_{n}\in \mathcal{H}_{n}^{1}$ then, as will be shown, $%
u_{n}$ can be chosen such that $u_{n}\left( x\right) \neq 0$. Under this
choice, the eigenvalue $\lambda _{n}$ can be computed from the equation
\begin{equation*}
\lambda _{n}\left \langle u_{n}\left( x\right) ,u_{n}\left( x\right) \right
\rangle =\left \langle u_{n}\left( x\right) ,Tu_{n}\left( x\right) \right
\rangle \text{.}
\end{equation*}

The following proposition gives an explicit formula for $\left \langle
u_{n}\left( x\right) ,Tu_{n}\left( x\right) \right \rangle $. But first, we
need to introduce one additional terminology.

Fix a unit vector $l_{0}\in S\left( x^{\bot }\right) \subset V$ and let $%
A_{l_{0}}\in Lie\left( T_{l_{0}}\right) $ be a vector such that the morphism
$\exp :[0,2\pi )\rightarrow T_{l_{0}}$ given by $\exp \left( \theta \right)
=e^{\theta A_{l_{0}}}$ is an isomorphism.

\begin{proposition}
\label{formula_prop}The following equation holds:%
\begin{equation}
\lambda _{n}\left \langle u_{n}\left( x\right) ,u_{n}\left( x\right) \right
\rangle =\int \limits_{0}^{\pi }\mu \left( \theta \right) \left \langle
u_{n}\left( x\right) ,e^{-\theta A_{l_{0}}}u_{n}\left( e^{\theta
A_{l_{0}}}x\right) \right \rangle d\theta ,  \label{formula_eq}
\end{equation}%
where $\mu \left( \theta \right) =\sin \left( \theta \right) /2$.
\end{proposition}

For a proof, see Appendix \ref{proofs_sec}$.$

Our strategy is to construct "good" section $u_{n}\in \mathcal{H}_{n}^{1}$
and then to use Equation (\ref{formula_eq}).

\subsubsection{Construction of a "good" section}

For every $n\geq 0$, choose a highest weight vector $\psi _{n}\in
\mathbb{C}
\mathcal{F}_{n}$ with respect to $\left( H,E,F\right) $ ( $H\psi _{n}=2n\psi
_{n}$). In addition, choose a highest weight vector $v_{1}\in
\mathbb{C}
V$ ($Hv_{1}=2v_{1}$).

For every $n\geq 1$, the section $\psi _{n-1}\otimes v_{1}$ is a highest
weight vector in $%
\mathbb{C}
\mathcal{V}_{n}$. In order to get a section in $\mathcal{H}_{n}^{1}$, first
step, apply the lowering operator $F$ and consider the section $\widetilde{u}%
_{n}=F^{n-1}\left( \psi _{n-1}\otimes v_{1}\right) \in \mathcal{V}_{n}^{1}$.

Let us denote by $P_{n-1}$ the weight zero spherical function $F^{n-1}\psi
_{n-1}\in \mathcal{F}_{n-1}^{0}$ and note that, under appropriate choice of
coordinates, $P_{n-1}$ is the classical spherical function on the sphere $%
S\left( V\right) $. The following proposition gives an explicit expression
of $\widetilde{u}_{n}$ in terms of the function $P_{n-1}$ and the vector $%
v_{1}$.

\begin{proposition}
\label{vector_prop}The section $\widetilde{u}_{n}$ can be written as%
\begin{equation}
\widetilde{u}_{n}=P_{n-1}\otimes v_{1}+\frac{1}{n}EP_{n-1}\otimes Fv_{1}+%
\frac{1}{2n\left( n+1\right) }E^{2}P_{n-1}\otimes F^{2}v_{1}\text{.}
\label{vector_eq}
\end{equation}
\end{proposition}

For a proof, see Appendix \ref{proofs_sec}$.$

The second step in order to get a section in $\mathcal{H}_{n}^{1}$ is to
define
\begin{equation*}
u_{n}\left( y\right) =P_{y}\widetilde{u}_{n}\left( y\right) ,
\end{equation*}%
for every $y\in S\left( V\right) $, where $P_{y}$ is the orthogonal
projector on $y^{\bot }$.

\begin{proposition}
\label{vector1_prop}The section $u_{n}$ is symmetric or anti-symmetric
depending on the parity of $n$ as follows%
\begin{eqnarray*}
u_{n} &\in &\mathcal{H}_{n}^{+,1}\text{ when }n\text{ is even,} \\
u_{n} &\in &\mathcal{H}_{n}^{-,1}\text{ when }n\text{ is odd.}
\end{eqnarray*}
\end{proposition}

For a proof, see Appendix \ref{proofs_sec}$.$

Now, we are ready to finish the computation. Using Equation (\ref{formula_eq}%
), we can write%
\begin{equation*}
\lambda _{n}\left \langle u_{n}\left( x\right) ,u_{n}\left( x\right) \right
\rangle =\int \limits_{0}^{\pi }\mu \left( \theta \right) \left \langle
e^{\theta A_{l_{0}}}u_{n}\left( x\right) ,\widetilde{u}_{n}\left( e^{\theta
A_{l_{0}}}x\right) \right \rangle d\theta .
\end{equation*}

Considering formula (\ref{vector_eq}), it is evident that the functions $%
EP_{n-1}$ and $E^{2}P_{n-1}$ must vanish at $x\in S\left( V\right) $ since
these are function of weight different then zero with respect to the action
of $T_{x}$. Hence
\begin{equation*}
u_{n}\left( x\right) =\widetilde{u}_{n}\left( x\right) =P_{n-1}\left(
x\right) v_{1}\in
\mathbb{C}
x^{\bot ,1}\text{.}
\end{equation*}

For $k=0,1,2$, define the integrals
\begin{eqnarray*}
I_{n}^{k} &=&\frac{1}{\left \langle u_{n}\left( x\right) ,u_{n}\left(
x\right) \right \rangle }\int \limits_{0}^{\pi }\mu \left( \theta \right)
\left \langle e^{\theta A_{l_{0}}}u_{n}\left( x\right) ,E^{k}P_{n-1}\left(
e^{\theta A_{l_{0}}}x\right) F^{k}v_{1}\right \rangle d\theta \\
&=&\frac{1}{P_{n-1}\left( x\right) \left \Vert v_{1}\right \Vert ^{2}}\int
\limits_{0}^{\pi }\mu \left( \theta \right) E^{k}P_{n-1}\left( e^{\theta
A_{l_{0}}}x\right) \left \langle e^{\theta A_{l_{0}}}v_{1},F^{k}v_{1}\right
\rangle d\theta \text{,}
\end{eqnarray*}%
and note that the eigen value $\lambda _{n}$ can be expressed as
\begin{equation}
\lambda _{n}=I_{n}^{0}+\frac{1}{n}I_{n}^{1}+\frac{1}{2n\left( n+1\right) }%
I_{n}^{2}\text{.}  \label{eigen_eq}
\end{equation}

\begin{theorem}[Main technical statement]
\label{technical_thm} For $k=0,1,2$, the integrals $I_{n}^{k}$ are equal to%
\begin{eqnarray*}
I_{n}^{0} &=&\left \{
\begin{array}{cc}
1 & n=1 \\
1/6 & n=2 \\
0 & n\geq 3%
\end{array}%
\right. , \\
I_{n}^{1} &=&\left \{
\begin{array}{cc}
0 & n=1 \\
-2/3 & n=2 \\
0 & n\geq 3%
\end{array}%
\right. , \\
I_{n}^{2} &=&\left \{
\begin{array}{cc}
0 & n=1 \\
0 & n=2 \\
2\left( -1\right) ^{n-1} & n\geq 3%
\end{array}%
\right. .
\end{eqnarray*}
\end{theorem}

For a proof, see Subsection \ref{tech_proof_sub}.

Consequently, using Theorem \ref{technical_thm} and Equation (\ref{eigen_eq}%
) we obtain the desired formula%
\begin{equation*}
\lambda _{n}=\frac{\left( -1\right) ^{n-1}}{n\left( n+1\right) }\text{,}
\end{equation*}%
which proves Theorem \ref{spectral_thm}.

\subsection{Proof of the main technical statement\label{tech_proof_sub}}

Let $\left( e_{1},e_{2},e_{3}\right) $ be an orthonormal basis of $V$.
Silently, the reader should think of the basis vector $e_{3}$ as standing in
place of the fixed unit vector $x\in S\left( V\right) $ and of the basis
vector $e_{2}$ as standing in place of the vector $l_{0}\in S\left( x^{\bot
}\right) $. It is possible to choose vectors $A_{e_{i}}\in Lie\left(
T_{e_{i}}\right) $ which satisfy the relations%
\begin{eqnarray*}
\left[ A_{e_{3}},A_{e_{1}}\right] &=&A_{e_{2}}, \\
\left[ A_{e_{3}},A_{e_{2}}\right] &=&-A_{e_{1}}, \\
\left[ A_{e_{1}},A_{e_{2}}\right] &=&A_{e_{3}},
\end{eqnarray*}%
and, in addition, satisfy $\left[ A_{e_{i}},A_{e_{j}}\right] =A_{e_{i}}e_{j}$
for every $1\leq i,j\leq 3$. We can define the following $sl_{2}$ triple $%
\left( H,E,F\right) $ which is associated with $T_{e_{3}}$

\begin{eqnarray*}
H &=&-2iA_{e_{3}}, \\
E &=&iA_{e_{2}}-A_{e_{1}}, \\
F &=&A_{e_{1}}+iA_{e_{2}}\text{.}
\end{eqnarray*}

\subsubsection{Spherical coordinates}

We introduce spherical coordinates $f:\left( 0,2\pi \right) \times \left(
0,\pi \right) \rightarrow S\left( V\right) $ given by $f\left( \varphi
,\theta \right) =g_{\varphi }\cdot \left( \cos \left( \theta \right)
e_{3}+\sin \left( \theta \right) e_{1}\right) $ where%
\begin{equation*}
g_{\varphi }=%
\begin{pmatrix}
\cos \left( \varphi \right) & -\sin \left( \varphi \right) & 0 \\
\sin \left( \varphi \right) & \cos \left( \varphi \right) & 0 \\
0 & 0 & 1%
\end{pmatrix}%
.
\end{equation*}

In the coordinates $\left( \varphi ,\theta \right) $ the operators $H,E,F$
are given by the following formulas%
\begin{eqnarray*}
H &=&2i\partial _{\varphi }, \\
E &=&-e^{-i\varphi }\left( i\partial _{\theta }+\cot \left( \theta \right)
\partial _{\varphi }\right) , \\
F &=&-e^{i\varphi }\left( i\partial _{\theta }-\cot \left( \theta \right)
\partial _{\varphi }\right) .
\end{eqnarray*}

\subsubsection{Highest weight vector in $V$}

The vector $v_{1}=-e_{1}+ie_{2}$ is a highest weight vector in $V$ and we
note that $\left \Vert v_{1}\right \Vert ^{2}=2$. For $k=0,1,2$, let us
denote by $j^{k}\left( \theta \right) $ the function $\left \langle
e^{\theta A_{e_{2}}}v_{1},F^{k}v_{1}\right \rangle $. Explicit calculation
reveals that%
\begin{eqnarray*}
j^{0}\left( \theta \right) &=&\cos \left( \theta \right) +1, \\
j^{1}\left( \theta \right) &=&2i\sin \left( \theta \right) , \\
j^{2}\left( \theta \right) &=&2\cos \left( \theta \right) -2.
\end{eqnarray*}

\subsubsection{Spherical function in $\mathcal{F}_{n}$}

For $n\geq 0$, let $P_{n}\in \mathcal{F}_{n}^{0}$ denote the unique weight
zero spherical function which satisfies the normalization condition $%
P_{n}\left( e_{3}\right) =1$. Define the generating function
\begin{equation*}
G\left( \varphi ,\theta ,t\right) =\sum \limits_{n=0}^{\infty }P_{n}\left(
\varphi ,\theta \right) t^{n}\text{.}
\end{equation*}

The generating function $G$ admits an explicit formula.

\begin{theorem}[\protect \cite{T}]
\begin{equation}
G\left( \varphi ,\theta ,t\right) =\left( 1-2t\cos \left( \theta \right)
+t^{2}\right) ^{-1/2}\text{.}  \label{generating_eq}
\end{equation}
\end{theorem}

\begin{remark}
Note that $G\left( 0,0,t\right) =\left( 1-t\right) ^{-1}=\sum
\limits_{n=0}^{\infty }t^{n}$ which is compatible with the normalization
condition $P_{n}\left( 0,0\right) =1$.
\end{remark}

Applying the raising operator $E$ we obtain the generating functions%
\begin{eqnarray*}
EG\left( \varphi ,\theta ,t\right) &=&\sum \limits_{n=0}^{\infty
}EP_{n}\left( \varphi ,\theta \right) t^{n}, \\
E^{2}G\left( \varphi ,\theta ,t\right) &=&\sum \limits_{n=0}^{\infty
}E^{2}P_{n}\left( \varphi ,\theta \right) t^{n}.
\end{eqnarray*}
Granting formula (\ref{generating_eq}), explicit calculation reveals that
\begin{eqnarray*}
EG\left( \varphi ,\theta ,t\right) &=&ie^{-i\varphi }t\sin \left( \theta
\right) \left( 1-2t\cos \left( \theta \right) +t^{2}\right) ^{-3/2}, \\
E^{2}G\left( \varphi ,\theta ,t\right) &=&-3e^{-2i\varphi }t^{2}\sin
^{2}\left( \theta \right) \left( 1-2t\cos \left( \theta \right)
+t^{2}\right) ^{-5/2}.
\end{eqnarray*}

\subsubsection{Putting everything together}

In terms of our choice of the highest weight vector $v_{1}$ and the
spherical functions $P_{n}$'s, the integrals $I_{n}^{k}$, $k=0,1,2$, are
given in the spherical coordinates $\left( \varphi ,\theta \right) $ by
\begin{eqnarray*}
I_{n}^{k} &=&\frac{1}{P_{n-1}\left( x\right) \left \Vert v_{1}\right \Vert
^{2}}\int \limits_{0}^{\pi }\mu \left( \theta \right) E^{k}P_{n-1}\left(
0,\theta \right) j^{k}\left( \theta \right) d\theta \\
&=&\frac{1}{2}\int \limits_{0}^{\pi }\mu \left( \theta \right)
E^{k}P_{n-1}\left( 0,\theta \right) j^{k}\left( \theta \right) d\theta .
\end{eqnarray*}

For $k=0,1,2$, define the generating functions
\begin{equation*}
I^{k}\left( t\right) =\sum \limits_{n=0}^{\infty }I_{n+1}^{k}t^{n}.
\end{equation*}

Each $I^{k}\left( t\right) $ can be expressed as the integral
\begin{equation*}
I^{k}\left( t\right) =\frac{1}{2}\int \limits_{0}^{\pi }\mu \left( \theta
\right) E^{k}G\left( 0,\theta ,t\right) j^{k}\left( \theta \right) d\theta .
\end{equation*}

Explicit calculation of the integrals $I^{k}\left( t\right) $ reveals that
\begin{eqnarray*}
I^{0}\left( t\right) &=&\frac{1}{2}\left( 1+\frac{1}{3}t\right) , \\
I^{1}\left( t\right) &=&\frac{1}{2}\left( -\frac{4}{3}t\right) , \\
I^{2}\left( t\right) &=&\frac{1}{2}\left( 4\left( 1+t\right)
^{-1}-4t-4\right) =2\sum \limits_{n=2}^{\infty }\left( -1\right) ^{n}t^{n}%
\text{.}
\end{eqnarray*}

From the above formulas, using Equation (\ref{eigen_eq}),we get
\begin{equation*}
\lambda _{n}=\frac{\left( -1\right) ^{n-1}}{n\left( n+1\right) },
\end{equation*}%
for every $n\geq 1$. This finish the proof of Theorem \ref{technical_thm}.

\appendix

\section{Proofs\label{proofs_sec}}

\subsection{Proof of Theorem \protect \ref{char_thm}}

Since $\tau $ is a morphism of $O\left( V\right) $ representations, it maps $%
V$ isometrically onto $\mathcal{H}_{1}^{-}$ - the unique anti-symmetric copy
( $\theta $ acts by $-1$) of the three dimensional representation of $%
SO\left( V\right) $, which, by Proposition \ref{dec_prop}, coincides with $%
\mathbb{V}=\mathcal{H}\left( \lambda _{\max }\right) $.

Evidently, $\tau $ is an isometry, up to a scalar. Hence, it is enough to
show that $Tr\left( \tau \circ \tau ^{T}\right) =3$, which we verify as
follows:%
\begin{eqnarray*}
Tr\left( \tau \circ \tau ^{T}\right) &=&r^{2}\cdot Tr\left( \alpha
_{can}\circ \alpha _{can}^{T}\right) =\frac{3}{2}\int \limits_{x\in S\left(
V\right) }Tr\left( i_{x}^{T}\circ i_{x}\right) dx \\
&=&\frac{3}{2}\int \limits_{x\in S\left( V\right) }2dx=3.
\end{eqnarray*}

Finally, the relation $\tau \circ i_{x}=\varphi _{x}$ follows from%
\begin{equation*}
ex_{x|\mathbb{V}}\circ \alpha _{can}=\mathrm{Pr}_{x}.
\end{equation*}

This concludes the proof of the theorem.

\subsection{Proof of Theorem \protect \ref{dec_thm}}

Property 1 is the classical result of spherical harmonics on the two
dimensional sphere, which can be found for example in \cite{T}. We just note
that the representation $\mathcal{F}_{n}$ consists of the restriction to $%
S\left( V\right) $ of harmonic polynomials of degree $n$, which implies that
$\mathcal{F}_{n}=\mathcal{F}_{n}^{+}$ when $n$ is even and $\mathcal{F}_{n}=%
\mathcal{F}_{n}^{-}$ when $n$ is odd.

Property 2 follows from property 1 since $\mathcal{N}$ can be trivialized
using the $O\left( V\right) $ invariant section $s\in \mathcal{N}$ where $%
s\left( y\right) =y$, for every $y\in S\left( V\right) $.

We now prove Properties 3 and 4 simultaneously.

Since $\mathcal{V}=\mathcal{F}\otimes V$ as a representation of $O\left(
V\right) $, we can compute the isotypic components of $\mathcal{V}$ in terms
of the isotypic components of $\mathcal{F}$. The computation proceeds as
follows:

\begin{itemize}
\item For $n=0$, $\mathcal{F}_{0}\otimes V=V^{-}$.

\item For $n\geq 1$, $\mathcal{F}_{n}\otimes V=\left( \mathcal{F}_{n}\otimes
V\right) _{n-1}^{?}\oplus \left( \mathcal{F}_{n}\otimes V\right)
_{n}^{?}\oplus \left( \mathcal{F}_{n}\otimes V\right) _{n+1}^{?}$ where
\begin{equation*}
?=\left \{
\begin{array}{cc}
+ & n\text{ odd} \\
- & n\text{ even}%
\end{array}%
\right. .
\end{equation*}

The decomposition of $\mathcal{F}_{n}\otimes V$ as a representation of $%
SO\left( V\right) $ is computed using the branching rules of a tensor
product and the action of $\theta \in O\left( V\right) $ is derived from the
facts that $V=V^{-}$ and Property 1.
\end{itemize}

This implies that the isotypic components of $\mathcal{V}$ are

\begin{itemize}
\item For $n=0$, $\mathcal{V}_{n}=\mathbf{1}^{+}$.

\item For odd $n\geq 1$, $\mathcal{V}_{n}=\left( \mathcal{F}_{n-1}\otimes
V\right) _{n}^{-}\oplus \left( \mathcal{F}_{n}\otimes V\right)
_{n}^{+}\oplus \left( \mathcal{F}_{n+1}\otimes V\right) _{n}^{-}$.

\item For even $n\geq 1$, $\mathcal{V}_{n}=\left( \mathcal{F}_{n-1}\otimes
V\right) _{n}^{+}\oplus \left( \mathcal{F}_{n}\otimes V\right)
_{n}^{-}\oplus \left( \mathcal{F}_{n+1}\otimes V\right) _{n}^{+}$.
\end{itemize}

Combined with Property 2 and the fact that $\mathcal{V=H}\oplus \mathcal{N}$
yields Properties 3,4.

This concludes the proof of the Theorem.

\subsection{Proof of Proposition \protect \ref{dec_prop}}

Fix $n\geq 1$. Denote $\mathcal{H=H}_{n}$ and $\mathcal{H}^{\pm }=\mathcal{H}%
_{n}^{\pm }$. The statement that $\mathcal{H}^{+}\subset \ker C$ follows
from the facts that $C\left( x,-y\right) =C\left( -x,y\right) =C\left(
x,y\right) $ and that a section $s\in $ $\mathcal{H}^{+}$ satisfies that $%
s\left( -x\right) =-\theta \left( s\right) \left( x\right) =-s\left(
x\right) $. Similarly, the statement that $\mathcal{H}^{-}\subset \ker O$
follows from the facts that $O\left( x,-y\right) =O\left( -x,y\right)
=-O\left( x,y\right) $ and that a section $s\in $ $\mathcal{H}^{-}$
satisfies that $s\left( -x\right) =-\theta \left( s\right) \left( x\right)
=s\left( x\right) $.

Since, by definition, $T=C-O$, this implies that $C_{|\mathcal{H}^{-}}=T_{|%
\mathcal{H}^{-}}$ and $-O_{|\mathcal{H}^{+}}=T_{|\mathcal{H}^{+}}$. Since $T$
commutes with the action of $SO\left( V\right) $ and $\mathcal{H}^{\pm }$
are irreducible representations
\begin{equation*}
T_{|\mathcal{H}^{\pm }}=\lambda ^{\pm }Id.
\end{equation*}

We are left to show that $\lambda ^{+}=\lambda ^{-}$. The argument proceeds
as follows:

Let us denote by $ev_{x}:%
\mathbb{C}
\mathcal{H\rightarrow
\mathbb{C}
}x^{\bot }$ the evaluation map at the point $x$. Since $x$ is fixed by the
group $T_{x}$, this implies that $ev_{x}$ is a morphism of representations
of $T_{x}$. Moreover, $ev_{x}$ induces an isomorphism of weight spaces%
\begin{equation*}
ev_{x}:\mathcal{H}^{\pm ,1}\overset{\simeq }{\longrightarrow }x^{\bot ,1}%
\text{.}
\end{equation*}

Fix sections $u^{\pm }$ $\in \mathcal{H}^{\pm ,1}$ normalized to be of norm
1. Applying Formula (\ref{formula_eq}), we can write%
\begin{eqnarray}
\lambda ^{\pm }\left \langle ev_{x}u^{\pm },ev_{x}u^{\pm }\right \rangle
&=&\int \limits_{0}^{\pi }\mu \left( \theta \right) \left \langle
ev_{x}u^{\pm },ev_{x}\left( e^{\theta A_{l_{0}}}u^{\pm }\right) \right
\rangle d\theta  \label{basic_eq} \\
&=&\left \langle ev_{x}u^{\pm },ev_{x}\left( \pi ^{\pm }\left( \overline{\mu
}\right) u_{n}\right) \right \rangle ,  \notag
\end{eqnarray}%
where $\pi ^{\pm }:T_{x}\rightarrow U\left(
\mathbb{C}
\mathcal{H}^{\pm }\right) $ are the group actions restricted to the subgroup
$T_{x}$ and $\overline{\mu }$ is the function on $T_{x}$ corresponding to $%
\mu $ via the isomorphism $e^{\theta A_{l_{0}}}$.

Equation (\ref{basic_eq}), implies that
\begin{equation*}
\lambda ^{\pm }=\left \langle u^{\pm },\pi ^{\pm }\left( \overline{\mu }%
\right) u^{\pm }\right \rangle _{%
\mathbb{C}
\mathcal{H}^{\pm }}.
\end{equation*}

This implies that $\lambda ^{\pm }$ are characterized solely in terms of the
irreducible representation $\pi ^{\pm }:SO\left( V\right) \rightarrow
U\left(
\mathbb{C}
\mathcal{H}^{\pm }\right) $, which, in turns, implies that $\lambda
^{+}=\lambda ^{-}$.

This concludes the proof of the proposition.

\subsection{Proof of Proposition \protect \ref{formula_prop}}

Let $f:T_{x}\times \left( 0,\pi \right) \rightarrow S\left( V\right) $ be
the spherical coordinates on $S\left( V\right) $ given by $f\left( g,\theta
\right) =ge^{\theta A_{l_{0}}}x$. In these coordinates, the normalized Haar
measure on $S\left( V\right) $ is given by $dg\boxtimes \mu \left( \theta
\right) d\theta $ where $dg$ is the normalized Haar measure on $T_{x}$ and $%
\mu \left( \theta \right) =\sin \left( \theta \right) /2.$

The section $u_{n}\mathcal{\in H}_{n}^{1}$ is a character vector with
respect to the group $T_{x}$, let us denote this character by $\chi
:T_{x}\rightarrow S^{1}$ and note that we have $g\cdot u_{n}=\chi \left(
g\right) u_{n}$, for every $g\in T_{x}$. Now, compute

\begin{eqnarray*}
\lambda _{n}\left \langle u_{n}\left( x\right) ,u_{n}\left( x\right) \right
\rangle &=&\left \langle u_{n}\left( x\right) ,Tu_{n}\left( x\right) \right
\rangle \\
&=&\int \limits_{y\in S\left( V\right) }\left \langle u_{n}\left( x\right)
,T\left( x,y\right) u_{n}\left( y\right) \right \rangle dy \\
&=&\int \limits_{T_{x}}dg\int \limits_{0}^{\pi }\mu \left( \theta \right)
\left \langle u_{n}\left( x\right) ,T\left( x,ge^{\theta A_{l_{0}}}x\right)
u_{n}\left( ge^{\theta A_{l_{0}}}x\right) \right \rangle d\theta \\
&=&\int \limits_{T_{x}}dg\int \limits_{0}^{\pi }\mu \left( \theta \right)
\left \langle u_{n}\left( x\right) ,gT\left( x,e^{\theta A_{l_{0}}}x\right)
g^{-1}u_{n}\left( ge^{\theta A_{l_{0}}}x\right) \right \rangle d\theta \\
&=&\int \limits_{T_{x}}dg\int \limits_{0}^{\pi }\mu \left( \theta \right)
\left \langle g^{-1}u_{n}\left( x\right) ,T\left( x,e^{\theta
A_{l_{0}}}x\right) g^{-1}u_{n}\left( ge^{\theta A_{l_{0}}}x\right) \right
\rangle d\theta \\
&=&\int \limits_{T_{x}}dg\int \limits_{0}^{\pi }\mu \left( \theta \right)
\left \langle u_{n}\left( x\right) ,T\left( x,e^{\theta A_{l_{0}}}x\right)
u_{n}\left( e^{\theta A_{l_{0}}}x\right) \right \rangle d\theta \\
&=&\int \limits_{0}^{\pi }\mu \left( \theta \right) \left \langle
u_{n}\left( x\right) ,e^{-\theta A_{l_{0}}}u_{n}\left( e^{\theta
A_{l_{0}}}x\right) \right \rangle d\theta ,
\end{eqnarray*}%
where, step 4 follows from the fact that $T$ commutes with the action of $%
SO\left( V\right) $ which is equivalent to the property that $T\left(
gx,gy\right) =gT\left( x,y\right) g^{-1}$, for every $x,y\in S\left(
V\right) $ and $g\in SO\left( V\right) $ which, in particular, implies that $%
T\left( x,ge^{\theta A_{l_{0}}}x\right) =T\left( gx,ge^{\theta
A_{l_{0}}}x\right) =gT\left( x,e^{\theta A_{l_{0}}}x\right) g^{-1}$ and step
7 follows from the fact that $T\left( x,e^{\theta A_{l_{0}}}x\right) $ is
the operator of parallel translation along the big circle connecting the
point $e^{\theta A_{l_{0}}}x$ with the point $x$.

This concludes the proof of the proposition.

\subsection{Proof of Proposition \protect \ref{vector_prop}}

First we note the following simple fact: The operator $EF:%
\mathbb{C}
\mathcal{V}_{n}\rightarrow
\mathbb{C}
\mathcal{V}_{n}$ preserve the weight spaces $\mathcal{V}_{n}^{l}$, and,
moreover, since $%
\mathbb{C}
\mathcal{V}_{n}$ is a representation of highest weight $2n$ with respect to
the $sl_{2}$ triple $\left( H,E,F\right) $ we have
\begin{equation}
EF_{|\mathcal{V}_{n}^{l}}=\left( n+l\right) \left( n-l+1\right) Id,
\label{weight_eq}
\end{equation}%
for $l=-n,..,n$.

Now, calculate%
\begin{equation*}
\widetilde{u}_{n}=F^{n-1}\left( \psi _{n-1}\otimes v_{1}\right) =\sum
\limits_{i=0}^{n-1}%
\begin{pmatrix}
n-1 \\
i%
\end{pmatrix}%
F^{n-1-i}\otimes F^{i}\left( \psi _{n-1}\otimes v_{1}\right) .
\end{equation*}

Since $%
\mathbb{C}
V$ is a representation of highest weight $2$ with respect to the $sl_{2}$
triple $\left( H,E,F\right) $, all tensors of the form $\left( -\right)
\otimes F^{k}v_{1}$, for $k\geq 3,$ vanish. This implies that the above sum
is equal to
\begin{equation*}
F^{n-1}\psi _{n-1}\otimes v_{1}+\left( n-1\right) F^{n-2}\psi _{n-1}\otimes
Fv_{1}+\frac{\left( n-1\right) \left( n-2\right) }{2}F^{n-3}\psi
_{n-1}\otimes F^{2}v_{1}.
\end{equation*}%
.

Recall that $P_{n-1}=F^{n-1}\psi _{n-1}$. Explicit calculation, using
formula (\ref{weight_eq}), reveals that
\begin{eqnarray*}
F^{n-1}\psi _{n-1} &=&\frac{1}{n\left( n-1\right) }EP_{n-1}, \\
F^{n-3}\psi _{n-1} &=&\frac{1}{\left( n-2\right) \left( n-1\right) n\left(
n+1\right) }E^{2}P_{n-1}.
\end{eqnarray*}

Combining all the above yields the desired formula for $\widetilde{u}_{n}$.

This concludes the proof of the proposition.

\subsection{Proof of Proposition \protect \ref{vector1_prop}}

The statement follows directly from the facts that $V=V^{-}$ which implies
that $\theta \left( F^{k}v_{1}\right) =-F^{k}v_{1}$ and that $P_{n-1}\in
\mathcal{F}_{n-1}^{?}$ where

\begin{equation*}
?=\left \{
\begin{array}{cc}
+ & n\text{ odd} \\
- & n\text{ even}%
\end{array}%
\right. .
\end{equation*}

This concludes the proof of the proposition.

\bigskip

\end{document}